\documentclass[12pt]{article}
\usepackage{amssymb}
\usepackage{amsmath}
\usepackage{mathtools}
\usepackage{amsthm}
\usepackage{amsbsy}
\usepackage{graphicx}
\usepackage[]{epstopdf}
\usepackage{tikz}
\usetikzlibrary{positioning}
\newcommand{\be}{\begin{equation}}
\newcommand{\ee}{\end{equation}}
\newcommand{\bea}{\begin{eqnarray}}
\newcommand{\eea}{\end{eqnarray}}
\newcommand{\ba}{\begin{array}}
\newcommand{\ea}{\end{array}}

\newcommand{\bc}{\begin{center}}
\newcommand{\ec}{\end{center}}
\newcommand{\ben}{\begin{enumerate}}
\newcommand{\een}{\end{enumerate}}
\newcommand{\bfi}{\begin{figure}}
\newcommand{\efi}{\end{figure}}

\newcommand{\bq}{\begin{quote}}
\newcommand{\eq}{\end{quote}}
\newcommand{\bqu}{\begin{quotation}}
\newcommand{\equ}{\end{quotation}}
\newenvironment{emphit}{\begin{itemize}}{\end{itemize}}
\newcommand{\bemp}{\begin{emphit}}
\newcommand{\eemp}{\end{emphit}}

\newcommand{\bt}{\begin{tabular}}
\newcommand{\et}{\end{tabular}}

\newtheorem{myth}{Theorem}[section]

\begin{document}
\date{}
\title{A MATHEMATICAL MODEL OF THE EFFECTS OF STRIKE ON NIGERIAN UNIVERSITIES
\footnote{2020 mathematics subject classification primary 34D05; secondary 34G10}
\thanks{{\bf keywords:  A mathematical model, Strike, Public Universities, Private Universities}}}
\author{A. O. Isere\thanks{All correspondence to be addressed to this author.}\\
Department of Mathematics,\\
Ambrose Alli University,\\
Ekpoma 310001, Nigeria.\\
isere.abed@gmail.com,\\
isereao@aauekpoma.edu.ng}\maketitle
\begin{abstract}
This paper carries out a scientific study of the effects of strike on Nigerian Universities. A mathematical model is formulated to examine the behavior of Nigerian University System when Public Universities are on strike. The results show that if all State and Federal Universities are on strike with the exception of the Private Universities, the University System in Nigeria is locally asymptotically stable.
\end{abstract}
\section{Introduction}
\paragraph{}
A mathematical model is an abstract description of a concrete system using mathematical concepts and language. Mathematical models can take many forms, differential equations, stochastic or game theoretic models. Mathematical models involving differential equations can be linear or non-linear and usually describes a system by a set of equations that establish relationship between the variables (see \cite{skd11, skd15, skd14}). Mathematical models find applications in Statistics, Computer Science, Physics, Biology, Earth Science, Chemistry and Engineering, generally, in physical and allied sciences.
Aside these traditional areas, mathematical models are used in social sciences and humanities. The use of mathematical models to solve problems in business or military operations is a large part of the field of operation research. Today mathematical models are used in epidemiology \cite{skd15,skd16}, music \cite{skd19}, linguistics and analytic philosophy \cite{skd17}. Also of interest is the use of mathematical models in the resolution of conflicts and economic cooperation \cite{skd41, skd12, skd13}.
\par{}
One of the earliest mathematical models is the malthusian growth model, sometimes called the simple exponential growth model. It stipulates that population grows exponentially or linearly. It is named after Thomas Robert Malthus who wrote an essay on the principle of population in $1798$.
\par{}
Malthusian models have the following forms:

$$ P(t)=P_{o}e^{rt} $$ which is the same as: $$ \frac{dP}{dt}=rP $$ with $$P(0)=P_{o}.$$
\par{}
In this paper, we will be using mathematical model in examining the dynamics of the university system in Nigeria. It is employing population models where the public universities are further subdivided into two subpopulations or compartments, the Federal ($U_{f}$), and the State Universities ($U_{s}$). The Private Universities ($U_{p}$) are considered as the third compartment of the student population in Nigerian University System. Therefore, this paper is carry out a scientific study of the dynamics (movement) of students from one compartment to the other, occasioned by strike, that is often not noticeable. For example, the Premium times \cite{skd30} stated that for the fifth year running (i.e up to $2020$), the university of Ilorin (UNILORIN) ranked as Nigeria's most preferred University by admission seekers in the country. But after about two decades, UNILORIN academic staff joined their colleagues in industrial action embarked upon by academic staff in March 2020 which lasted for months. From the year $2021$ admission record released by JAMB, as of August, 2022, they could only have 19 percent of 13,634 admission in $2020$ \cite{skd30}. Students prefer to go to universities where their education will not be interrupted by strike. B.O Monogbe and T. G. Monogbe in 2019  carried out an empirical investigation on ASUU strike and Nigerian educational system. The study examines the extent to which strikes influences tertiary educational system in Nigeria. Their study further reveals that the quality of education and the students performance is negatively influenced by the incessant strike \cite{skd21}. Punch of August 14, 2020 reports that Varsities suspend 461,745 students' admission over strike, and further reports that a total of 36,947 candidates of the 2021 Unified Tertiary Matriculation Examination processed their admission into private institutions.  Private investment in education is not new, Kitaev in $2003$ asserted that private investment in education has become a popular policy in developed countries, namely, France, Britain, Germany, Spain, Singapore, Russia, Canada, Australia, Italy, China Japan, etc. In Africa, private universities are embraced to serve an as alternative to public universities in the region \cite{skd18}. This is where the problem lies. It should not be seen as an alternative but rather complimentary. But, in Nigeria, due to population explosion and incessant industrial action embarked upon by the staff of public universities, people therefore are advocating for the inclusion of private universities in the country, stated by Ssewamela in 2014 as reported in \cite{skd18}. Several studies have been conducted on private universities in Nigeria, those studies dwelled on the emergence and the analysis of Private Universities as stated in Obasi 2007.
\newpage
\begin{table}[!hbp]
\begin{center}
\begin{tabular}{|c||c|c|c|}
\hline
S/N & YEAR & PRIVATE UNIVERSITIES & STRIKE DURATION IN DAYS  \\
\hline
1 & 1999-2000 & 3 & 150  \\
\hline
2 & 2001-2002 & 4 & 104  \\
\hline
3 & 2003-2004 & 1 & 180 \\
\hline
4 & 2005-2006 & 15 & 21 \\
\hline
5 & 2007-2008 & 10 & 97 \\
\hline
6 & 2009-2010 & 7 & 157 \\
\hline
7 & 2011-2012 & 9 & 59 \\
\hline
8 & 2013-2014 & 0 & 165 \\
\hline
9 & 2015-2016 & 19 & 0 \\
\hline
10 & 2017-2018 & 7 & 120 \\
\hline
11 & 2019-2020 & 4 & 270 \\
\hline
12 & 2021-2022 & 32 & 210 \\
\hline
          & {\bf TOTAL} & {\bf 111} & {\bf 1,323} \\
\hline
\end{tabular}
\end{center}
\caption{The Table showing Private universities and Strike Duration}\label{quatable4}
\end{table}

To mention that research to date has shown that most of the studies conducted on Private Universities are done quantitatively \cite{skd18}. Thus, there is need to qualitatively study the university system in Nigeria, therefore, this work will carry out a qualitative study of the effect of strike on Nigerian University System using a mathematical model.
\section{The Formulation of the Mathematical Model}
\begin{equation}
\frac{dU_f}{dt}=\Lambda_f +\alpha_{sf}U_s +\alpha_{pf}U_p -dU_f -\alpha_{fs}U_f - \alpha_{fp}U_f -\lambda_f U_f
\end{equation}
\begin{equation}
\frac{dU_s}{dt}=\Lambda_s +\alpha_{ps}U_p +\alpha_{fs}U_f -dU_s -\alpha_{sf}U_s - \alpha_{sp}U_s -\lambda_s U_s
\end{equation}
\begin{equation}
\frac{dU_p}{dt}=\Lambda_p +\alpha_{sp}U_s +\alpha_{fp}U_f -dU_p -\alpha_{ps}U_p - \alpha_{pf}U_p -\lambda_p U_p
\end{equation}
$$ U_f (0)= N, U_{s}(0)=0, U_P (0)=0 $$
Where

  N = total human population  eligible for University for University education.\\

 $ U_{f} = $ total population of students in federal Universities.\\

 $ U_s = $ total population of students in state Universities.\\

 $ U_p = $ total population of students in private Universities.\\

  $\Lambda_f = $ rate of admission of students into federal Universities.\\

 $ \Lambda_s = $ rate of admission of students into state Universities.\\

 $ \Lambda_p = $ rate of admission of students into private Universities.\\

 $ \alpha_{sf} = $ rate of movement of students from state to federal Universities\\

 $\alpha_{sf} = $ rate of movement of students from state to federal Universities\\

 $\alpha_{pf} = $ rate of movement of students from private to federal Universities largely due to finance at time t\\

 $\alpha_{fs} = $ rate of movement of students from federal to some state Universities largely due to strike at time t\\

 $\alpha_{fp} = $ rate of movement of students from federal to private Universities largely due to strike at time t\\

 $\alpha_{ps} = $ rate of movement of students from private to state Universities largely due to finance at time t\\

 $\alpha_{sp} = $ rate of movement of students from some state to private Universities largely due to strike at time t\\

 $ d = $ the natural death rate at time t\\

 $\lambda_f = $ rate of graduation of students from federal Universities at time t.\\

 $\lambda_s = $ rate of graduation of students from state Universities at time t.\\

 $\lambda_p = $ rate of graduation of students from private Universities at time t\\

 \par{} The quantity $\alpha_{sf}U_s$ represents the number of students who move from state to federal universities at a given time $t$ and $\alpha_{fs}U_f$ represents the number of students who move from federal to state universities at a given time $t$ due to finance and strike respectively. Also, $\alpha_{fp}U_f$ and $\alpha_{pf}U_p$ represent the number of students that move from federal to private and from private to federal universities largely due to ASUU strike and finance respectively.
Similarly, $\alpha_{sp}U_s$ and $\alpha_{ps}U_p$ represent the number of students who move from state to private and from private to state universities respectively. Each of them is a vector quantity with the subscripts indicating the direction of movement.

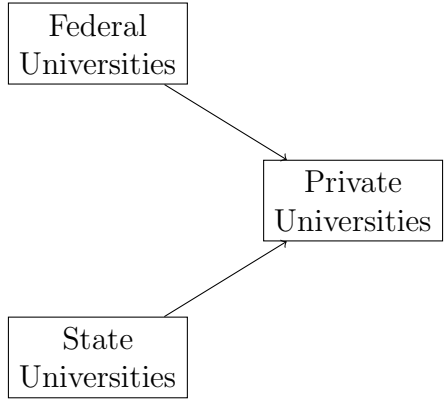
\begin{figure}[!htb]
\begin{center}
\begin{tikzpicture}[every node/.style=
{rectangle,draw=black,align=center}]\node(c){Private\\ Universities};
\node(a)[above left =of c]{Federal\\ Universities} edge [->](c);
\node(b)[below left =of c]{State\\ Universities} edge [->](c);
\end{tikzpicture}
\caption{A Schematic Diagram For the Three Compartments}\label{skd001}
\end{center}
\end{figure}

\subsection{EQUILIBRIUM STATE}
Equilibrium analysis gives the fixed points or equilibrium solutions. That is, the values of the university population for which the system will no longer experience movement. At that equilibrium point all the rate of movement are equated to zero. However, during strike, the rate of productivity or graduation of student is dormant or equated to zero for $U_f$ and for most $U_s$ but other rates are very active. And these are of interest to us. Interestingly also, is the particular time $t$ when there were no students in the universities, when the university system was at the formative years. This happened in 1948 when the first University was founded.  We refer to such state as the zero state.
$$ \frac{dU_f}{dt}= N-d $$
$$\frac{dU_s}{dt}= 0 $$ and $$ \frac{dU_p}{dt}= 0$$
or the movement free equilibrium.\\
Setting the RHS of equation (1)-(3) to zero and solving algebraically, we have an equilibrium point $(U^{*}_{f},U^*_{s},U^*_{p})$, where
\begin{equation}
U^*_{f}= \frac{\Lambda_f + \alpha_{sf}U^*_s + \alpha_{pf}U^*_{p} }{d + \alpha_{fs} + \alpha_{fp} + \lambda_f}
\end{equation}
\begin{equation}
U^*_{s}= \frac{\Lambda_s + \alpha_{fs}U^*_f + \alpha_{ps}U^*_{p} }{d + \alpha_{sf} + \alpha_{sp} + \lambda_s}
\end{equation}
\begin{equation}
U^*_{p}= \frac{\Lambda_p + \alpha_{fp}U^*_f + \alpha_{sp}U^*_{s} }{d + \alpha_{pf} + \alpha_{ps} + \lambda_p}
\end{equation}

\subsection{STABILITY ANALYSIS}
Local stability analysis helps to determine the behaviour of the different compartments near the equilibrium solution.
We do this by obtaining the Jacobian matrix of the system.
\par{}
The Jacobian of the above equations (1)-(3) is:
$$ J =
\begin{bmatrix}
-d -\alpha_{fs} -\alpha_{fp} -\lambda_{f} &  \alpha_{sf} &   -\alpha_{pf}\\
 \alpha_{fs} & -d -\alpha_{sf} -\alpha_{sp} -\lambda_{s} & -\alpha_{ps}\\
 \alpha_{fp} &  \alpha_{sp}&  -d -\alpha_{ps} -\alpha_{pf} -\lambda_{p}
\end{bmatrix}
$$
During strike, it is assumed that there is no movement of students from private universities to state or federal universities and that no students move from state to federal universities. Consequently, $\lambda_{f}=0$ and also $\lambda_{s}=0 $ except a situation where some state universities are partially or fully in session during strike, then $\lambda_{s}$ and $\lambda_{p}$ are not zeros. The stability analysis of the various situations will be presented below.
\begin{myth}
Given that all state universities are partially or fully in session together with all private universities. Then, the Nigerian University System is locally asymptotically stable.
\end{myth}
{\bf Proof:}\\
Since all state universities are partially or fully in session during ASUU strike. That is $\alpha_{ps}=0$, $\alpha_{pf}=0$ and $\alpha_{sf}=0$.
$$ J =
\begin{bmatrix}
-d -\alpha_{fs}-\alpha_{fp} &  0 &   0\\
 \alpha_{fs} & -d -\lambda_{s} & -\alpha_{ps}\\
 \alpha_{fp} &  0 &  -d -\alpha_{ps} -\lambda_{p}
\end{bmatrix}
$$
The eigenvalues ($\lambda_{1},\lambda_{2},\lambda_{3}$) are: $\lambda_{1}=-d-\lambda_{s}$, $\lambda_{2}=-d-\alpha_{ps}-\lambda_{p}$ and $\lambda_{3}=-d-\alpha_{fs}-\alpha_{fp}$
Thus, the Nigeria University System is locally asymptotically stable.

\begin{myth}
Given that some state universities are partially in session during ASUU strike, then, Nigerian university system is
locally asymptotically stable.
\end{myth}
{\bf Proof:}\\
Since some state universities are partially or fully in session during ASUU strike. That is $\alpha_{ps}=0$, $\alpha_{pf}=0$ and $\alpha_{sf}=0$.
The Jacobian matrix of the proposed model in (1)-(3) evaluated at ASUU strike period gives:

$$ J_{0} =
\begin{bmatrix}
-d -\alpha_{fs} &  0 &   0\\
 \alpha_{fs} & -d -\alpha_{sp} -\lambda_{s} & 0\\
 \alpha_{fp} &  0 &  -d  -\lambda_{p}
\end{bmatrix}
$$

Eigenvalues of $J_{0}$ obtained from the characteristic polynomial \begin{equation} \lambda^3 -a_{1}\lambda^{2} -a_{1}\lambda + a_{3}=0  \end{equation} where
\begin{gather*}
a_{1} =-3d-\lambda_{p}-\alpha_{sp}-\lambda_{s}-\alpha_{fs}\\
a_{2} =-3d^{2}-2d\alpha_{sp}-2d\lambda_{p}-2d\lambda_{s}-\alpha_{fs}\alpha_{sp}\\
         -\alpha_{fs}\lambda_{p}-\alpha_{fs}\lambda_{p}-\alpha_{fs}\lambda_{s}-\alpha_{sp}\lambda_{p}-\lambda_{p}\lambda_{s}\\
a_{3} = (d+\lambda_{p})(d+\alpha_{sp}+\lambda_{s})(d+\alpha_{fs})
\end{gather*}
are:\\
$\lambda_{1}=-d-\lambda_{p}$, $\lambda_{2}=-d-\alpha_{sp}-\lambda_{s}$ and $\lambda_{3}=-d-\alpha_{fs}$ \\
 Hence eigenvalues ($\lambda_{1},\lambda_{2},\lambda_{3}$) of $J_{0}$ are negative. Therefore, Nigerian University System is locally asymptotically stable provided some (all) state and private universities are in section in Nigeria during ASUU strike.

\begin{myth}
Suppose that all state and federal universities are on strike with the exception of the private universities in Nigeria.
The university system in Nigeria is locally asymptotically stable.
\end{myth}
{\bf Proof:}\\
Setting other rates to zero except $\alpha_{sp}$ and $\alpha_{fp}$; and also $\lambda_{f}$ and $\lambda_{s}$ to zero, the jacobian matrix becomes:

$$ J(p) =
\begin{bmatrix}
-d  &  0 &   0\\
 0 & -d -\alpha_{sp} & 0\\
 \alpha_{fp} &  0 &  -d  -\lambda_{p}
\end{bmatrix}
$$
The eigenvalues of the characteristic polynomial \begin{equation} \lambda^3 -a_{1}\lambda^{2} -a_{2}\lambda + a_{3}d=0  \end{equation}

where:
$a_{1}=3d-\lambda p -\alpha_{sp}$, $ a_{2}= -3d^{2}-2d\alpha_{sp}-2d\lambda_{p}-\alpha_{sp}\lambda_{p}$ and
$a_{3}=(d+\lambda_{p})(d+\alpha_{sp})$ are: $\lambda_{1}=-d$, $\lambda_{2}=-d-\lambda_{p}$ and $\lambda_{3}=-d-\alpha_{sp}$. \\ Hence, the eigenvalues of $J(p)$ are all negative. Thus, the Nigerian University System is still locally asymptotically stable.
\par{}
This accounts for the reason while ASUU strike can linger over time, since the presence of private universities gives Nigerian University System an asymptotic stability.

\subsection{The Movement Restricted State Equilibrium}
The movement restricted state is the equilibrium $(\frac{U_{f}}{d}, 0, 0)$. Trivially, this state is always stable
and there will be no movement from one university compartment to another during strike. Then, ASUU strike will be efficacious.
\par{} This is what happened prior to 1999 when there was no private universities, and once ASUU  was on strike, all public universities complied and Government had no option than to attend to ASUU very quickly.
\par{}
In that state, the expected population of students in the universities would be:
\begin{equation}
\frac{dU_{f}}{dt}= \Lambda_{f}-dU_{f}
\end{equation}
This is equivalent to the following:
\begin{equation}\label{sd1}
\frac{dU}{dt}= \Lambda-dU
\end{equation}
since there was no private universities. So, solution of (\ref{sd1}) gives:
\begin{equation}
U(t)= \Lambda/d  + (U_{i}-\Lambda/d)e^{-dt}
\end{equation}
where $U_{i}$ is the initial number of students in public universities in Nigeria. It is easy to see that as
$t\rightarrow \infty$ then $U\rightarrow \Lambda/d$, which is the asymptotic population size of university students at the movement restricted state. Hence, the entire population of University students will be in one University System.
\par{}
{\bf Heart burning questions:} With the advent of state universities in Nigeria since 1982 and private universities in 1999, will there still be movement restricted state? If your answer is No, then what is the implication of that on ASUU strike?
\section{Concluding Remarks}
So long as movement within the university compartments during strike exits, the University System in Nigeria will continue to be asymptotically stable and strike will yield little or no intended results.
\section{Recommendation}
Public universities should avoid a situation that can push $\lambda_{f}$ and $\lambda_{s}$ to zero which usually happens during strike. In other words, $\alpha_{fp}$ and $\alpha_{sp}$ particularly should be made dormant or redundant by ensuring that the public universities are on.

\end{document}